# Differences between real and integer production possibility sets in data envelopment analysis

**Dariush Khezrimotlagh**

Department of Applied Statistics, Faculty of Economics and Administration, University of Malaya, Kuala Lumpur, Malaysia, dariush@um.edu.my, January 20, 2015.

## Abstract

The real Production Possibility Set (PPS) is completely generated by observations in the real approach, but generating the integer PPS may not be possible by only using observations in the integer approach. In other words, this phenomenon says that all points in the real generated PPS are dominated by a point of the linear combination of observations, whereas there might be some points in the integer generated PPS which are not dominated with the points of the linear combination of observations. This paper shows how the integer production possibility set is made and the mixed-integer linear programming is defined. The paper also addresses some shortcomings in the recent mixed-integer linear programming while the integer axioms are considered.

## 1. Introduction

Khezrimotlagh *et al.* [1] noted some of the shortcomings on the integer-valued radial model in Data Envelopment Analysis (DEA), and supported their claims with a counterexample in [2], about two mixed integer redial models called Lozano and Villa Model (LVM) [3] and Kuosmanen and Kazemi-Matin Model (KKM) [4-5]. However, recently a chapter book [6] has been written which demonstrates the same shortcomings. The following sections illustrate the Mathematical drawbacks of the inaccurate discussions and Theorem 1 in [4-6].

## 2. Can a non-integer variable be used in an integer axiom?

One of the proposed axiom in DEA is the convexity axiom. The real convexity axiom says that "if two points A and B of a set T (with real values) are selected, then all the points of the line segment AB (with real values) belong to T", that is, if $(x, y) \in T \subseteq \mathbb{R}_+^{m+p}$ and $(x', y') \in T \subseteq \mathbb{R}_+^{m+p}$, then $[\lambda(x, y) + (1 - \lambda)(x', y')] \in T \subseteq \mathbb{R}_+^{m+p}$, for $\lambda \in [0,1]$.

Therefore, the integer convexity axiom can be defined as "if two points A and B of a set T (with integer values) are considered, then all the points of the line segment AB (with integer values) belong to T", that is, if $(x, y) \in T \subseteq \mathbb{Z}_+^{m+p}$ and $(x', y') \in T \subseteq \mathbb{Z}_+^{m+p}$, then $[\frac{u}{v}(x, y) + (1 - \frac{u}{v})(x', y')] \in T \subseteq \mathbb{Z}_+^{m+p}$, for $u, v \in \mathbb{Z}_+$, where $u \leq v$, $v|x - x'$ and $v|y - y'$.

Replacing $\lambda$ with two integer values $u$ and $v$, is mistakenly interpreted as restricting $\lambda$ to the set of rational numbers in [6]. However, this is the divisibility of the set of integer numbers. In fact, $[\frac{u}{v}(x, y) + (1 - \frac{u}{v})(x', y')]$, that is, $\frac{u}{v}(x - x', y - y') + (x', y')]$ is integer while $v$ divides both $x - x'$ and $y - y'$. This condition, that is, $v$ divides both $x - x'$ and $y - y'$, is the properties of the divisibility of the set of integer numbers, and $\lambda = u/v$ should not be interpreted as restricting $\lambda$ to the set of rational numbers.

There should not be non-integer variables in an integer axiom and there are not any in the above integer convexity axiom. *What's the meaning of "integer convexity axiom" if one generates a non-integer value with a linear combination of integer values?* As the next section illustrates, using the above axiom is avoided in [4-6], which effects their improper further conclusions.

## 3. Is the formulation of KKM valid?

The following elementary example clearly depicts how a real (an integer) Production Possibility Set (PPS) is generated with the real (integer) DEA axioms and how Mathematical equations are defined to generate the DEA PPS, which is a rudimentary step to understanding the base of DEA.

Suppose there are two DMUs A(5, 9) and B(2, 2). Figure 1 (Figure 2) depicts using the first real (integer) axiom, that is, *inclusion of observation*. There is no difference between Figures 1 and 2, since the components of A and B are integer in both approaches. Indeed, the set of integer numbers is a subset of the set of real numbers, that is, $\mathbb{Z} \subseteq \mathbb{R}$.

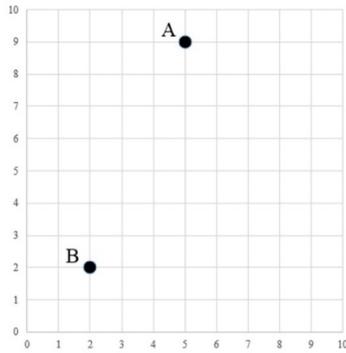 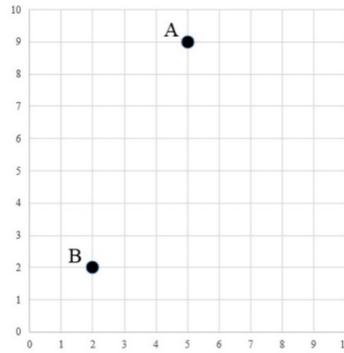

Figure 1: Using the first real DEA axiom.     Figure 2: Using the first integer DEA axiom.

Now, let's suppose that the second DEA axiom, *convexity axiom*, is applied. In other words, let's consider the linear combinations of A and B, that is, $\lambda(x_A, y_A) + (1 - \lambda)(x_B, y_B)$, for $\lambda \in [0,1]$, which can be written as $\lambda_A(x_A, y_A) + \lambda_B(x_B, y_B)$, where $\lambda_A + \lambda_B = 1$, and $\lambda_A \geq 0$ and $\lambda_B \geq 0$. It yields finding the points $(\lambda_A x_A + \lambda_B x_B, \lambda_A y_A + \lambda_B y_B)$, where $\lambda_A + \lambda_B = 1$, and $\lambda_A \geq 0$ and $\lambda_B \geq 0$.

As Figure 3 shows, the points of the line segment AB are completely generated by the real convexity axiom. However, there are no points with integer-valued components on the line segment AB as Figure 4 illustrates.

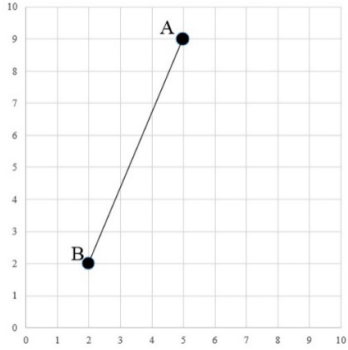
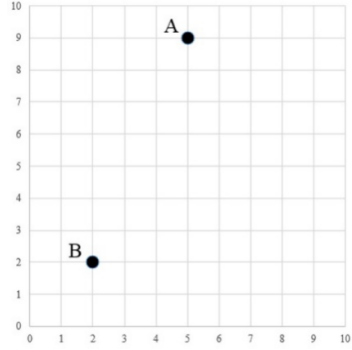

Figure 3: Using the second real DEA axiom.    Figure 4: Using the second integer DEA axiom.

Now, let's apply the third real (integer) axiom, free disposability, that is, if $(x, y) \in \mathbf{T}$, then $(x', y) \in \mathbf{T}$, for $x' \geq x$, and $(x, y') \in \mathbf{T}$, for $y' \leq y$, where $x, x', y$ and $y'$ are real (integer). In other words, let's consider the following linear inequalities where $\lambda_A + \lambda_B = 1$, $\lambda_A \geq 0$ and $\lambda_B \geq 0$ (indeed, combination of both axioms yields the following linear equations which are the base of formulations in DEA):

$$\lambda_A x_A + \lambda_B x_B \leq x', \qquad \lambda_A y_A + \lambda_B y_B \geq y'.$$

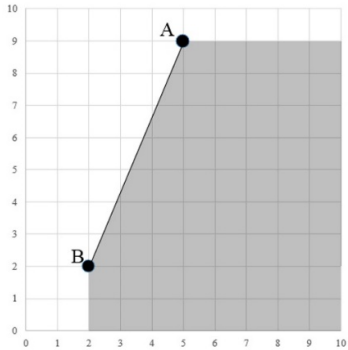
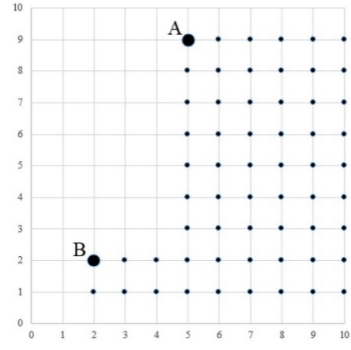

Figure 5: Using the third real DEA axiom.    Figure 6: Using the third integer DEA axiom.

As depicted in Figure 5, applying the third real axiom completes the PPS in the real approach, however, the PPS in the integer approach is not yet generated (Figures 6, 7).

One may argue, why isn't the third axiom (free disposability) used in the second step of generating the integer PPS (as it is usually defined as a second axiom in the literature of DEA or the axioms are simultaneously considered)? Indeed, after that in the third step, the second axiom (convexity axiom) can completely generate the integer PPS!

This argument is valid, however, the mathematical combination of such approaches is not linear and none of the LVM and KKM (nor any of DEA models) are formulated with this approach. In fact, for such an aim, the following non-linear equations should have been considered to generate the integer PPS:

$$x_1 \geq x_A, \qquad x_2 \geq x_B, \qquad y_1 \leq y_A, \qquad y_2 \leq y_B,$$
$$\lambda_A x_1 + \lambda_B x_2 = x', \qquad \lambda_A y_1 + \lambda_B y_2 = y',$$
$$x_1 \in \mathbb{Z}_+, \quad x_2 \in \mathbb{Z}_+, \quad x' \in \mathbb{Z}_+, \quad y_1 \in \mathbb{Z}_+, \quad y_2 \in \mathbb{Z}_+, \quad y' \in \mathbb{Z}_+,$$
$$\lambda_A + \lambda_B = 1, \lambda_A \geq 0 \text{ and } \lambda_B \geq 0.$$

*One may also wish that "the set of axioms assumed to generate the DEA technology apply simultaneously, not sequentially (one after another)" or "applying the axioms to any feasible activity not just to observed units".* However, none of the current mathematical DEA formulations support such wishes. As it is clearly illustrated "applying the axioms simultaneously" causes the formulation to be non-linear. Of course, any feasible activities should adapt to the applied axioms, but it does not mean that all feasible activities are considered in the linear combination of inputs and/or outputs in DEA formulation. It should not have been forgotten that in mathematical formulation of a DEA model, only the linear combination of observed DMUs are considered. Indeed, convexity axiom is only applied for observed values in a DEA model. If one hopes to consider the linear combination of all feasible activities, the mathematical formulation is not linear, it is exactly non-linear as the above formulation illustrates it!

These obvious statements are neglected in [4-6] while attempting to generate the integer PPS and formulate KKM according to the integer DEA axioms. Thus, the following lemma is proved (Figures 2, 4, 6-8) which rejects Theorem 1 in [4-6].

**Lemma**: Suppose that the integer convexity axiom is considered. Then, $\mathbf{T}_{\text{IDEA}}^{\text{VRS}}$, may not be equal to $\mathbf{T}_{\text{DEA}}^{\text{VRS}} \cap \mathbb{Z}_+^{m+p}$, where $(x_i, y_i) \in \mathbb{Z}_+^{m+p}$ are observed DMUs ($i = 1, 2, \dots, n$), $m$ number of inputs and $p$ number of outputs:

$$\mathbf{T}_{\text{IDEA}}^{\text{VRS}} = \left\{ (x, y) \in \mathbb{Z}_+^{m+p} : x \geq \sum_{i=1}^n \lambda_i x_i, y \leq \sum_{i=1}^n \lambda_i y_i, \sum_{i=1}^n \lambda_i = 1 \right\},$$

and

$$\mathbf{T}_{\text{DEA}}^{\text{VRS}} = \left\{ (x, y) \in \mathbb{R}_+^{m+p} : x \geq \sum_{i=1}^n \lambda_i x_i, y \leq \sum_{i=1}^n \lambda_i y_i, \sum_{i=1}^n \lambda_i = 1 \right\}.$$

**Corollary**:

$$\mathbf{T}_{\text{IDEA}}^{\text{VRS}} = \left\{ \begin{matrix} (x, y) \in \mathbb{Z}_+^{m+p} : x = \sum_{i=1}^n \lambda_i x_i', \ y = \sum_{i=1}^n \lambda_i y_i', x_i' \geq x_i \\ y_i' \leq y_i, (x_i', y_i') \in \mathbb{Z}_+^{m+p}, \sum_{i=1}^n \lambda_i = 1, \ \lambda_i \geq 0 \end{matrix} \right\} = \mathbf{T}_{\text{DEA}}^{\text{VRS}} \cap \mathbb{Z}_+^{m+p}.$$

In order to generate the remaining points shown in Figure 7, it is enough to select the non-observed point E and then apply the second integer axiom to generate points C and D (Figure 8). After that, applying the third integer axiom completes the integer PPS (Figure 9). Indeed, the following constraints yield generation of the remaining points shown in Figure 8, where $\lambda_B + \lambda_E = 1, \lambda_B \geq 0, \lambda_E \geq 0, x' \in \mathbb{Z}_+$ and $y' \in \mathbb{Z}_+$:

$$\lambda_B x_B + \lambda_E x_E \leq x', \qquad \lambda_B y_B + \lambda_E y_E \geq y'.$$

The above indisputable statements demonstrate that the real PPS is completely generated by observations in the real approach, but generating the integer PPS may not be possible by only using observations in the integer approach. In other words, this phenomenon says that all points in the real generated PPS are dominated by a point of the linear combination of observations, whereas there might be some points in the integer generated PPS which are not dominated with the points of the linear combination of observations (Figure 7).

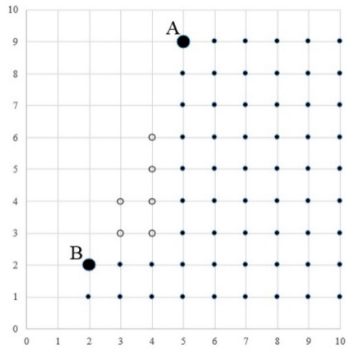

Figure 7: Non-generated points.

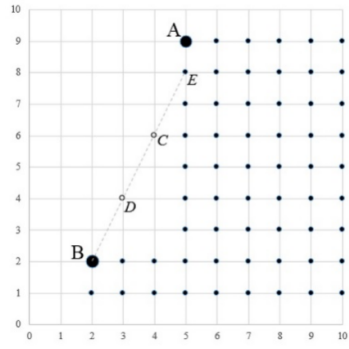

Figure 8: Using the third integer DEA axiom.

Now, let's suppose that $F$ **(5, 6)** is also observed (Figure 9). The following linear combination of observations clearly generates the point $C$ **(4, 6)** where $\frac{4}{9} + \frac{1}{3} + \frac{2}{9} = \mathbf{1}$:

$$\frac{4}{9}x_A + \frac{1}{3}x_B + \frac{2}{9}x_F = \mathbf{4}, \qquad \frac{4}{9}y_A + \frac{1}{3}y_B + \frac{2}{9}y_F = \mathbf{6}.$$

It means that generating the point $C$ **(4, 6)** is possible by observations A, B and F according to the integer convexity axiom, whereas it is impossible to generate $C$ **(4, 6)** by observations A and B (excluding $F$) as shown in Figures 6-8. Indeed, generating $C$ **(4, 6)** by observations A and B only, means that using the real convexity axiom instead using the integer convexity axiom (Figures 3-6).

Therefore, if the integer convexity axiom is accepted, it must also be accepted that $\lambda_F \neq \mathbf{0}$, which changes the meaning of optimal intensity weights.

The above clear arguments neither depend on whether the intensity weights $\lambda_A$, $\lambda_B$ and $\lambda_F$ are unique, nor depend on the use of returns to scale technologies. If the integer convexity axiom (or other axioms corresponding to the returns to scale) is used in KKM, its results are equivalent with the results of LVM (as it is claimed in [1]). Indeed, the simple outcome of the convexity axiom is that every valid linear combination of elements of a PPS must belong to the PPS.

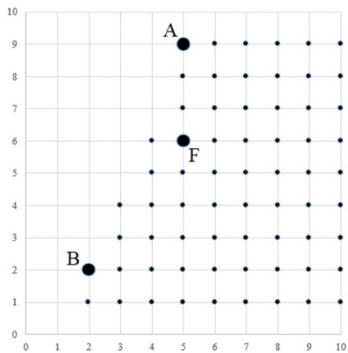

Figure 9: The generated integer PPS.

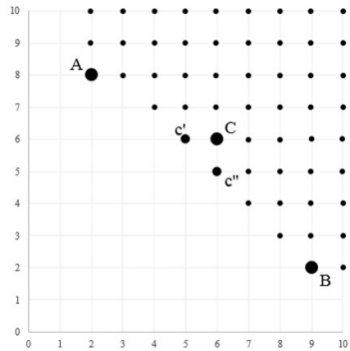

Figure 10: Over-estimating efficiency.

KKM generates the real PPS and after that it looks for integer values, which is completely against the attempt of generating the integer PPS according to the introduced integer axioms. In such an approach, its results should not be compared as advantages via

LVM. As a result, the formulation of KKM is not valid and the following theorem is proved.

**Theorem**: Regarding the integer axioms, the formulation of KKM is not valid.

Note that if one claims that the generated integer PPS can have non-integer values, there is no logical reason to propose the integer axioms and such PPS should not be called an integer PPS. In fact, no other boundaries can be considered or generated while the input and output values are restricted to the set of integer numbers. As it is argued in [1] there is no need to introduce new axioms while generating the real PPS and its intersection with integer values are considered. *If such an approach is valid to have a lower efficiency score, why does one not use CCR? Because rounding targets in this case has fewer disadvantages instead of collapsing the integer convexity axiom.*

On the other hand, it is argued in [4-6] that, the first equality constraint in LVM relaxes the free disposability axiom, which is correct in the integer approach only. However, that equality guarantees the integer convexity axiom and all points which should have been made by the integer free disposability axiom are dominated with observations (in other words, none of observations benchmark to their dominated points). Therefore, relaxing this axiom does not have effect on the LVM optimization. In contrast, replacing inequality by equality generates the real PPS and does not produce the integer PPS as discussed above, which is completely against the integer axioms and effects the KKM optimization.

## 4. Does KKM over-estimate efficiency?

Although, the LVM formulation over-estimates efficiency, the KKM formulation, even with generating the invalid PPS also over-estimates efficiency.

In order to clarify this point, let's suppose DMUs A(2, 8, 1), B(9, 2, 1) and C(6, 6, 1) as depicted in Figure 10. Both LVM and KKM technical efficiency scores are 1 for C, however, it is clear that both models over-estimate efficiency. Indeed, C is an inefficient DMU which is dominated by $c'$ and $c''$. In other words, even by expanding an inappropriate PPS, KKM still over-estimates efficiency.

The problem of over-estimating efficiency is a problem of all redial approaches which are neither able to simultaneously measure the potential decreasing of inputs and increasing of outputs, nor able to benchmark DMUs toward the economical part of the efficient frontier [7-8].

Nonetheless, even accepting the invalid real PPS by KKM in the integer approach, the counterexample in [5] still rejects this claim that "the KKM input's targets are not greater than the LVM input's targets". Although, in [6] KKM attempts to explain that "the optimal intensity weights and slacks are not necessarily unique", the robust counterexample in [5] does not depend on this simple and clear statement. In fact, the optimal intensity weights yield that the linear combination of first inputs to be $\sum_{i=1}^{5} \lambda_i^* x_{i1} = 37.5$, which would never allow KKM to suggest 37 as LVM does.

Note that, the point (37, 234; 250, 100) in the counterexample in [5] belongs to the integer PPS. Indeed, the intensity weights $\lambda_1 = 0.05263158$, $\lambda_2 = 0.05263158$, $\lambda_3 = 0.05263158$, $\lambda_4 = \lambda_5 = 0$, yield that the linear combination of first inputs to be $\sum_{i=1}^{5} \lambda_i x_{i1} = 36.84210526$, which allows KKM to suggest 37, however, these intensity weights are not optimum. This is an obvious weakness of the redial approaches while PPS is expanded and can even be seen by finding the optimum input targets of the CCR model, which are $(37.5, 233.\overline{3})$ where $\lambda_3^* = 1.3\overline{8}$ and $\lambda_1^* = \lambda_2^* = \lambda_4^* = \lambda_5^* = 0$. In other words, expanding a PPS by a redial approach does not certainly yield finding a valid efficient target.

## 5.   Examples in [6]

There are two simple examples in [6] which are not correct unfortunately. In the first example two DMUs $A(5, 12; 3)$ and $B(10, 12; 2)$ are considered. Then, the third axiom is applied to find $A'(6,12; 3)$. Here is the gap which was illustrated in the above Section 3 and the remainder of the illustration in [6] is not valid. Indeed, KKM is not formulated with such approach, i.e., "using the free disposability before using the convexity axiom" or "using the free disposability and convexity axioms simultaneously".

In the second example, three DMUs $A(2,1)$, $B(3,2)$ and $C(3,1)$ are considered. Then, it is illustrated that the MILP algorithm will arbitrarily identify two different classes of slacks to benchmark C, that is, $(1,0)$ and $(0,1)$. However, none of DEA models and DEA algorithms select $(1,0)$ as appropriate slacks. Indeed, the optimum of lambdas never allows any of DEA algorithms to select $(1,0)$. It is not allowed to find the optimum slacks without considering the optimum lambdas. Therefore, this example is not valid clearly, and does not show a disadvantage of the additive measure.

## 6.   Conclusion

This note shows that the KKM formulation by accepting the integer axioms is not valid. Moreover, it rejects Theorem 1 proposed in [4-6].